\documentclass{amsart}
\usepackage{amssymb,amscd}
\usepackage[french]{babel}
\usepackage[latin1]{inputenc}
\usepackage[T1]{fontenc}
\newcounter{intro}

\newtheorem{theo}[intro]{Théorème}
\newtheorem{coro}[intro]{Corollaire}

\newtheorem{thm}{Théorème}[section]
\newtheorem{lem}[thm]{Lemme}
\newtheorem{prop}[thm]{Proposition}
\newtheorem{cor}[thm]{Corollaire}

\theoremstyle{remark}
\newtheorem{rem}[thm]{Remarque}
\newtheorem*{merci}{Remerciements}

\numberwithin{equation}{section}

\newcommand{\cref}[1]{Corollary~\ref{#1}}

\newcommand{\R}{\mathbb{R}}

\newcommand{\bS}{\mathbb{S}}

\newcommand{\cH}{\mathcal{H}}

\def\tq{{\rm\ \ tel\ que\ \ }}
\def\et{{\rm\ et\ }}
\def\si{{\rm\ si\ }}

\DeclareMathOperator{\ima}{Im}

\DeclareMathOperator{\rang}{rang}

\DeclareMathOperator{\vol}{vol}

\def\lra{\longrightarrow}

\begin{document}

\title[Cohomologie $L^2$ et parabolicité]
{Cohomologie $L^2$ et parabolicité}

\author{Gilles Carron}
\address{Laboratoire de Math\'ematiques Jean Leray (UMR 6629), 
Universit\'e de Nantes,
2, rue de la Houssini\`ere, B.P.~92208, 44322 Nantes Cedex~3, France}
\email{Gilles.Carron@math.univ-nantes.fr}

\begin{abstract}
Nous donnons ici une interprétation topologique des espaces de formes harmoniques $L^2$ de variétés
riemanniennes qui sont des déformations compactes d'espaces symétriques à courbure négative ou nulle
ou de groupes de Lie nilpotent simplement connexe. Nous étudions aussi la
cohomologie $L^2$ des variétés dont les bouts sont paraboliques, nous retrouvons en particulier des
résultats de M. Atiyah, V. Patodi, I. Singer, de W. Müller et une petite partie de résultats récents
de T. Hausel, E. Hunsicker, R. Mazzeo.
\end{abstract}
\maketitle
\section{Introduction}

Lorsque $(M,g)$ est une variété riemannienne complète, nous notons $\cH^k(M)$ l'espace des $k$-formes sur $M$ qui
sont $L^2$ et harmoniques :
$$\cH^k(M)=\{\alpha\in L^2(\Lambda^kT^k M), d\alpha=0 \et d^*\alpha=0\}.$$

Grâce aux cèlèbres théorèmes de Hodge et de-Rham, nous savons que pour une variété compacte, ces espaces
sont isomorphes aux groupes de cohomologie réel de la variété. Pour une variété non compacte, il est
naturel de se demander :
\begin{itemize}
\item[i)] A quelles conditions, les espaces de formes harmoniques $L^2$ sont de dimension finies ?
\item[ii)] Et quelle est alors l'interprétation topologique de ces espaces ?
\end{itemize}
Nous connaissons désormais de nombreuses réponses à ces questions (cf. par
exemple \cite{APS,carma,cargafa,HHM,M,MP,SS,Y,Z}).
Nous montrerons ici les deux résultats suivants :
\begin{theo} Soit $(N,g)$ un groupe de Lie nilpotent simplement connexe de dimension $n\ge 3$ et $(M,g)$
une variété riemannienne isométrique à l'infini
\footnote{On dit que deux variétés riemaniennes $(X_1,g_1)$ et $(X_2,g_2)$ sont isométriques à l'infini
s'il existe deux compacts $K_1\subset X_1$ et $K_2\subset X_2$ tels que $(X_1\setminus K_1, g_1)$ et
$(X_2\setminus K_2, g_2)$ sont isométriques.} à plusieurs copies de $(N,g)$ alors
$$\cH^k(M)\simeq \left\{\begin{array}{ll}
H^k_c(M)& {\rm si  \ } k \le n-2 \\
H^k(M)& {\rm si  \ } k \ge 2 \\ \end{array}
\right. .$$
\end{theo}
Ce théorème est une généralisation d'un résultat obtenu dans \cite{carma} à propos des variétés euclidiennes à
l'infini. Pour les degrés compris entre $2$ et $n-2$ la preuve de ce résultat est basée sur des
considérations élémentaires. Pour les degrés $1$ et $n-1$, nous devrons notamment utiliser un résultat de J. Jost et
K. Zuo (\cite{JZ}). Nous obtiendrons aussi avec une preuve similaire le résultat suivant :
\begin{theo} Soit $X=G/K$ un espace symétrique à courbure négative ou nulle (de dimension $n\ge 3$)
 et $(M,g)$ une variété riemannienne isométrique à l'infini
à plusieurs copies de X alors lorsque $\rang G\not =\rang K$ nous avons 
$$\cH^k(M)\simeq \left\{\begin{array}{ll}
H^k_c(M)& {\rm si  \ } k \le n-2 \\
H^k(M)& {\rm si  \ } k \ge 2 \\ \end{array}
\right. .$$
Lorsque $\rang G=\rang K$, alors cette conclusion persiste exceptée en degré $k=\dim X/2$, où $\cH^k(M)$
est de dimension infinie.
\end{theo}
 Remarquons que d'après un résultat de J. Lott (\cite{Lo}), la finitude des espaces de formes harmoniques $L^2$
 considérés dans ces deux théorèmes est une conséquence de la nullité de ces espaces pour les groupes de Lie
 nilpotents simplement connexe (\cite{caroma}) et pour les espaces symétriques (\cite{B}). J. Lott montre
  en effet que si deux variétés orientables $(X_1,g_1)$ et $(X_2,g_2)$ sont isométriques à l'infini alors 
$$  \dim \cH^k(X_1,g_1)<\infty \Longleftrightarrow \dim \cH^k(X_2,g_2)<\infty.$$
Nos hypothèses ne permettent pas de traiter le cas où $M$ est une surface riemannienne dont les bouts sont 
 euclidiens. En fait pour une surface riemannienne l'espace des $1$-formes harmoniques $L^2$ est
un invariant conforme et pour une telle surface $(S,g)$ de type fini nous avons : 
\begin{itemize}
 \item Soit $\cH^1(S,g)$ est de dimension finie et $(S,g)$ est conformément équivalente à une surface
 compacte privée d'un nombre fini de points et $\cH^1(S,g)\simeq \ima[H^1_c(S)\lra H^1(S)]$.
 \item Soit $\cH^1(S,g)$ est de dimension infinie.
 \end{itemize}
 La différence entre ces deux cas pouvant être aussi distingué par la parabolicité ou la non-parabolicité
  de la surface. Nous
 démontrerons le résultat partiellement analogue suivant :
\begin{theo} On considère $(M,g)$ une variété riemannienne qui au dehors d'un compact $D$ est isométrique à
$$( ]0,\infty[ \times\partial D,dt^2+h_t) $$
 où $h_t$ est une
 famille de métrique riemannienne sur $\partial D$ telle que : \begin{itemize}
 \item si $t>s$ alors $h_t\ge h_s$,
 \item si on note $L(t)=\vol_{n-1}( \{t\}\times \partial D)=\int_{\partial D}\sqrt{\det_{h_0} h_t}
 d\vol_{h_0}$ alors $$\int_0^\infty \frac{dt}{L(t)}=\infty ,$$
\item si on note  $J(t,\theta)=\sqrt{\det_{h_0(\theta)} h_t(\theta)}$ alors la fonction $$j(t)=\max_{\theta\in \partial D} J(t,\theta)\int_0^t \frac{1}{J(s,\theta)} ds$$ vérifie 
$$\int_0^\infty \frac{dt}{j(t)}=\infty.$$
\end{itemize} Alors nous avons l'isomorphisme :
$$\cH^k_2(M)\simeq \ima [H^k_c(M)\lra H^k(M)].$$
\end{theo}
Ce théorème général nous permet de donner une preuve unifiée de beaucoup de résultats connus par exemple
lorsque $(M,g)$ est à bouts cylindriques on retrouve un résultat de M.Atiyah, V.Patodi et I.Singer 
(\cite{APS}); nous retrouvons aussi un résultat de W. Müller à propos des $b$-métriques sur les variétés dont le bord
à un coin de codimension $2$ (\cite{Mu}); enfin nous re obtenons une partie d'un résultat récent de 
T.Hausel, E.Hunsicker et R.Mazzeo (\cite{HHM}) ; et enfin ce résultat permet de démontrer :
\begin{coro} Soit $E\lra X$ un fibré vectoriel riemannien sur une variété riemannienne compacte où les
fibres sont des plans euclidiens, $L$ équipé de sa métrique riemannienne vérifie :
 $$\cH^k_2(L)\simeq \ima [H^k_c(L)\lra H^k(L)].$$
\end{coro}
Ceci permet répond à une question posée dans \cite{caroma}.

\begin{merci} Je tiens à remercier l'Université de Stanford pour son hospitalité et particulièrement Rafe
Mazzeo dont les questions pertinentes sont à l'origine de ce travail. Je remercie aussi Colin Guillarmou et Nader Yeganefar pour
l'attention qu'ils ont porté à ce travail. Enfin ce travail a été fait alors que je bénéficiais d'une
délégation partielle au C.N.R.S.et d'une aide d'un programme ACI du ministère de la recherche.
\end{merci}

\section{Cohomologie $L^2$ et cohomologie usuelle}

L'objectif de ce paragraphe est d'utiliser des arguments élémentaires de
cohomologie usuelle pour donner une interprétation topologique des espaces de
formes harmoniques $L^2$ de certaines variétés riemanniennes (non compactes).
Nous commençons par rappeler brièvement les propriétés principales de la
cohomologie $L^2$ réduite\footnote{Désormais, on omettra de signaler qu'il s'agit de la cohomologie $L^2$
réduite.} ; on trouvera plus de détails dans l'article de J. Lott
\cite{Lo} ou encore dans \cite{caroma}.
Soit $(M^n,g)$ une variété riemannienne, on note $H^k_2(M)$ son $k^{ieme}$
espace de cohomologie $L^2$ , il est défini par
$$H^k_2(M)=\frac{\{\alpha\in L^2(\Lambda^kT^*M),
d\alpha=0\}}{\overline{\{d\beta,\ \beta\in L^2(\Lambda^{k-1}T^*M)\et d\beta\in
L^2(\Lambda^kT^*M)\}} }$$

Il est vrai que la cohomologie $L^2$ peut être calculée par des formes
$C^\infty$, i.e.  toute classe contient une forme lisse. Par exemple lorsque
$(M,g)$ est complète (ce que l'on supposera désormais) ces espaces s'identifient aux espaces de formes harmoniques
$L^2$.
$$H^k_2(M)\simeq \cH^k(M)=\{\alpha\in L^2(\Lambda^kT^*M), d\alpha=0 \et d^*\alpha=0\}.$$
En effet, si nous notons :
\begin{eqnarray*}
&Z^k_2(M)&=\{\alpha\in L^2(\Lambda^kT^*M), d\alpha=0\}\\
& &=\{\alpha\in L^2(\Lambda^kT^*M),
\tq \forall \beta\in C^\infty_0(\Lambda^{k+1}T^*M)\ \ \langle\alpha,d^*\beta\rangle=0\}\\
 \et\ &B^k_2(M)&=\overline{d C^\infty_0(\Lambda^{k-1}T^*M)}\\ \end{eqnarray*}
alors nous avons $H^k_2(M)\simeq \cH^k(M)=Z^k_2(M)\cap \left(B^k_2(M)\right)^\perp.$
Lorsque $\Omega\subset M$ est un ouvert à bord compact lisse de $M$, on peut définir aussi la cohomologie
$L^2$ absolue ou relative avec :
\begin{eqnarray*}&Z^k_2(\Omega)\! \! \! &=\{\alpha\in L^2(\Lambda^kT^*\Omega),
\tq \forall \beta\in C^\infty_0(\Lambda^{k+1}T^*\Omega)\ \ \langle\alpha,d^*\beta\rangle=0\}\\
\et \! \! \!&B^k_2(\Omega)&=\overline{d C^\infty_0(\Lambda^{k-1}T^*\overline{\Omega})}\\
\end{eqnarray*}
où une forme de $C^\infty_0(\Lambda^{k-1}T^*\overline{\Omega})$ est à support compact dans
l'adhérence de $\Omega$, elle n'est pas forcément nulle dans un voisinage de $\partial \Omega$.
Alors on définit $$H^k_2(\Omega)=Z^k_2(\Omega)/B^k_2(\Omega).$$ 
L'inclusion $j\,:\, \Omega\lra M$ induit donc
une application linéaire $$j^*\,:\, H^k_2(M)\lra H^k_2(\Omega).$$
La cohomologie $L^2$ relative est définie avec les espaces :
\begin{eqnarray*}
&Z^k_2(\Omega,\partial \Omega)&=\{\alpha\in L^2(\Lambda^kT^*\Omega),
\tq \forall \beta\in C^\infty_0(\Lambda^{k+1}T^*\overline{\Omega})\ \ \langle\alpha,d^*\beta\rangle=0\}\\
\et &B^k_2(\Omega,\partial \Omega)&=\overline{d C^\infty_0(\Lambda^{k-1}T^*\Omega)}\\
\end{eqnarray*}
et on définit $H^k_2(\Omega,\partial \Omega)=Z^k_2(\Omega,\partial \Omega)/B^k_2(\Omega,\partial \Omega)$.
L'application extension par zéro permet de définir une application linéaire
$e\,:\,  H^k_2(\Omega,\partial \Omega)\lra H^k_2(M)$.
Bien sûr si $\Omega$ est relativement compact ces espaces s'idenfient à la cohomologie absolue ou relative de $\Omega$.
Dans \cite{Lo} (voir aussi \cite{caroma}), il est remarqué que si $D$ est un domaine
compact de $M$ alors les suites exactes longues 
$$
..\lra H^k(D,\partial D)\stackrel{\mbox{e}}{\longrightarrow} H^k(M)\stackrel{\mbox{j}^*}{\longrightarrow} H^k(M\setminus D)
\stackrel{\mbox{b}}{\longrightarrow} H^{k+1}(D,\partial D)\lra..$$

$$..\lra H^k(M\setminus D,\partial D)\stackrel{\mbox{e}}{\longrightarrow} H^k(M)\stackrel{\mbox{j}^*}{\longrightarrow} H^k( D)
\stackrel{\mbox{b}}{\longrightarrow} H^{k+1}(M\setminus D,\partial D)\lra.. $$
 induisent deux petites suites exactes :
 $$H_2^k(M)\stackrel{\mbox{j}^*}{\longrightarrow} H_2^k(M\setminus D)
\stackrel{\mbox{b}}{\longrightarrow} H^{k+1}(D,\partial D)$$
 \begin{equation}
 \label{suiteex}
 H_2^k(M\setminus D,\partial D)\stackrel{\mbox{e}}{\longrightarrow}
 H_2^k(M)\stackrel{\mbox{j}^*}{\longrightarrow} H^k( D)
 \end{equation}
 
 Un autre héritage de ces suites exactes est le suivant :
 
 \begin{lem} 
 \label{inj1}
 Si $D$ est un domaine compact de $M$ et si $H^{k-1}(M\setminus
 D)=\{0\}$ alors
 $$\{0\}\rightarrow H^k(D,\partial D)\rightarrow H^k_2(M).$$
 \end{lem}
 
 \proof Soit $c\in H^k(D,\partial D)$, on sait que $c$ contient une forme fermée lisse
 à support compact dans l'interieur de $D$. Si $\alpha$ est nulle en cohomologie
 $L^2$, nous savons qu'elle est nulle en cohomologie usuelle
 (\cite{dR}, théorème 24). Ainsi $$c=[\alpha]\in \ker \left(\mbox{i}\,:\,H^k(D,\partial
 D)\longrightarrow H^k(M)\right)=\ima \left(\mbox{b}\,:\,H^{k-1}(M\setminus
 D)\longrightarrow H^k(D,\partial D)\right),$$ or ce dernier espace est nul par
 hypothèse.
 \endproof
 On a aussi compte tenu de (\ref{suiteex})
 \begin{lem} \label{inj2}
  Si $D$ est un domaine compact de $M$ et si $$H_2^{k}(M\setminus
 D,\partial D)=\{0\}$$ alors
 $$\{0\}\rightarrow H^k_2(M)\rightarrow H^k(D).$$
 \end{lem}
Remarquons que lorsque $H^{k}(\partial D)=\{0\}$ alors en degré $k$ la cohomologie relative
de $D_0$ se surjecte dans la cohomologie usuelle de $D_0$ et en combinant  les
deux lemmes précédents on obtient :
 \begin{prop}\label{iso}
 Supposons que
 pour un domaine compact $D\subset M$, on ait 
 $H_2^{k}(M\setminus D,\partial D)=\{0\}$ ainsi que $H^{k-1}(M\setminus
 D)= \{0\}$ et $H^{k}(\partial D)=\{0\}$ 
 alors on a
 $$H^k_2(M)\simeq H^k(D)\simeq H^{k}( D,\partial D).$$
 Et si de plus $H^{k}(M\setminus
 D)= \{0\}$ alors on a $H^k_2(M)\simeq H^k(M)\simeq H^{k}( D,\partial D)$.
 \end{prop}
 
Nous voulons maintenant un critère qui garantisse que
 $H_2^{k}(M\setminus D,\partial D)$ soit nul :
 \begin{lem} \label{nul}
 Supposons que $H^k_2(M)=\{0\}$ et que $H^k(\partial D)=\{0\}$ alors
 $$ H_2^{k}(M\setminus D)=\{0\}.$$ 
 \end{lem}
 \proof  Soit $c\in  H_2^{k}(M\setminus D)$, et $\alpha$ un représentant lisse de
 $c$, comme $H^k(\partial D)=\{0\}$, on peut trouver $\beta$ une $(k-1)$-forme lisse à
 support compact dans l'adhérence de $M\setminus D$ , de façon à ce que
$ \alpha-d\beta$, qui est encore un élément de $c$, soit à support compact dans
$M\setminus D$. Mais cette forme est nulle en cohomologie $L^2$ sur $M$,
on peut donc trouver une suite $(\varphi_l)_l$ de formes lisses à support compact dans $M$ tel
que 
$$\alpha-d\beta=L^2\lim_{l\to\infty} d\varphi_l.$$
En particulier en restriction à $M\setminus D$, on a 
$$\alpha-d\beta=L^2\lim_{l\to\infty} d(j^*\varphi_l).$$
donc $c$ est nulle.
 \endproof
 
 La concaténation de ces résultats nous donne le théorème suivant :
 \begin{thm} 
 \label{theogen}
 Soit $(M_0^n,g_0)$ une variété riemannienne complète orientée avec $D_0\subset
 M$ un domaine compact de $M$ tel que
  \begin{itemize}
  \item $H^{k-1}(M_0\setminus D_0) = H^{k}(M_0\setminus D_0)=\{0\}.$
  \item $H^{k-1}(\partial D_0)=H^{k}(\partial D_0)=\{0\}.$
  \item $H_2^{k}(M_0)=\{0\}.$
 \end{itemize}
 alors pour toute variété riemannienne $(M,g)$ isométrique au dehors d'un
 compact à $(M_0\setminus D_0,g_0)$, on a
 $$H^k_2(M)\simeq H^k(M)$$
 \end{thm}
 En effet , d'abord la dualité induite par l'opérateur de Hodge implique que 
 $H^{n-k}(\partial D_0)=\{0\}$ et $H_2^{n-k}(M_0)=\{0\}$; donc grâce au lemme (\ref{nul})
  on a $H_2^{n-k}(M_0\setminus D_0)=\{0\}$. La dualité induite par l'opérateur de Hodge nous permet
   d'affirmer $H_2^{n-k}(M_0\setminus D_0)\simeq H_2^{k}(M_0\setminus D_0,\partial D_0)$ et nous pouvons alors conclure grâce à la proposition (\ref{iso}).

 Une application immédiate est la suivante :
 \begin{cor} Soit $g_0$ une métrique complète sur $\R ^n$ et $k\in[2,n-2]$ un entier tel que $(\R^n,g_0)$ n'ait pas de 
 $k$-formes harmoniques $L^2$ non triviale  alors
pour toute variété riemannienne $(M,g)$ isométrique à $(\R^n,g_0)$ au dehors d'un compact on a:
$$H^k_2(M,g)\simeq H^k(M)\simeq  H_c^k(M).$$
 \label{coreu}\end{cor}
 \begin{rem}
 Lorsque $k$ vaut $0$ ou $n$, on a simplement $H^0_2(M,g)\simeq H_c^0(M)=\{0\}$ et
 $H^n_2(M,g)\simeq H^n(M)=\{0\}$.
 Le cas où $(M,g)$ est isométrique à l'infini à  plusieurs copies de $(\R^n,g_0)$ donne le même résultat.
  Nous verrons plus loin comment traité, avec un peu d'analyse, le cas des degré $1$ ou $n-1$.
 \end{rem}

Ce dernier corollaire a beaucoup d'applications :

\begin{itemize}
\item Si $g_0$ est la métrique euclidienne, on re-prouve un résultat de \cite{carma} ;
 d'ailleurs dans ce papier, on utilisait déjà l'injectivité (\ref{inj1}).
\item Si $(N,g)$ est un groupe de Lie nilpotent simplement connexe équipé d'une métrique
invariante à gauche alors $N$ est difféomorphe à $\R^n$ et on sait que $(N,g)$ n'a pas de
formes harmoniques $L^2$ non triviale (cf. corollaire 2.4 de \cite{caroma}).
\item Si $X=G/K$ est un espace symétrique à courbure négative ou nulle, 
alors $X$ est
difféomorphe à $\R^n$ et les travaux de A. Borel (\cite{B}) nous apprennent que
lorsque $\rang G\not=\rang K$ alors $X$ n'a pas de formes harmoniques $L^2$ non trivial, et que lorsque
 $\rang G=\rang K$ et $k\not=n/2$, $X$ n'a pas de $k$-formes harmoniques $L^2$ non
triviales.
\end{itemize}
\section{$L^2$ cohomologie et parabolicité}
\subsection{Un résultat d'annulation.}
Le but de cette partie est de prolonger les arguments de la partie précédente en y introduisant un
peu de géométrie. Nous commençons tout d'abord par le lemme suivant qui est un leger raffinement
dans une formulation un peu différente d'un résultat de
  P. Li, J. Jost et K. Zuo, J. Mc Neal et N. Hitchin \cite{Li,JZ,Mc,H}.
\begin{lem} 
\label{poid}
Soit $(M,g)$ une variété riemannienne complète et $\alpha\in L^2(\Lambda^pT^*M)$ une
$p$ forme $L^2$ fermée. On suppose qu'il existe $\beta\in C^\infty (\Lambda^{p-1}T^*M)$ avec
$$\alpha=d\beta\  \et\  \int_M \frac{|\beta|^2}{\psi(\rho)^2} \,d\vol <\infty $$  où
$\rho$ est la fonction distance à un point $o$ fixé de $M$ et 
$\psi$ est une fonction continue positive telle que
$$\int_1^\infty\frac{dr}{\psi(r)}=\infty$$ alors
la classe de $\alpha$ en cohomologie $L^2$ est nulle.
\end{lem}
\proof
Soient $r,R$ deux nombres réels tels que $0<r<R$, on leurs associe la fonction $\phi_{r,R}$ définie
par
$$\phi_{r,R}(s)=\left\{\begin{array}{lll}
1& \si & s\le r\\
\int_s^R \frac{dt}{\psi(t)} \times  \left(\int_r^R \frac{dt}{\psi(t)}\right)^{-1} &\si &s\in[r,R]\\
0& \si &s\ge R\\
\end{array}
\right.$$
La forme $\phi_{r,R} (\rho)\beta$ est lipschitsienne à support compact et on a
$$\alpha-d(\phi_{r,R}(\rho) \beta)=\alpha-\phi_{r,R}(\rho) d\beta- \phi'_{r,R}(\rho)d\rho \wedge
\beta$$
mais nous avons 
$$\|\phi'_{r,R}(\rho)d\rho \wedge
\beta \|^2_{L^2}\le  \left(\int_r^R \frac{dt}{\psi(t)}\right)^{-2} \int_{B(o,R)\setminus B(o,r)}
\frac{|\beta|^2}{\psi(\rho)^2} d\vol$$
 Or par hypothèse on  peut trouver deux suites $r_k<R_k$ telles que $\lim_{k\to\infty } r_k=\infty$ et 
 $$\left(\int_{r_k}^{R_k} \frac{dt}{\psi(t)}\right)^{-1}  \le \frac{1}{k}$$ d'où pour ces suites
 $$\alpha=L^2\lim_{k\to\infty} d (\phi_{r_k,R_k} (\rho)\beta).$$\endproof
 \begin{rem}
 \label{jost}
 Si on fait comme J. Jost et K. Zuo l'hypothèse qu'il existe une constante $C$ telle que
 $$\int_{B(o,R)}|\beta |^2d\vol\le C R^2$$ alors nos hypothèses sont satisfaites pour la fonction 
 $\psi(r)=r\log r$.
  \end{rem}
 Nos hypothèses sont intimement liées à la parabolicité. Pour plus d'information concernant cette
 notion, je conseille le survol remarquable et très complet de A. Grigor'yan \cite{gri}. Je
 rappelle ici juste que si $(M,g)$ est une variété riemannienne complète et $\mu$ une mesure
 absolument continue par rapport à $d\vol$ alors on dit que l'opérateur $\Delta^\mu$ associé
 à la forme quadratique $f\mapsto \int_M |df|^2 d\mu$ est parabolique \footnote{Si $d\mu=d\vol$,
 on dit simplement que $(M,g)$ est parabolique.} si l'une des conditions équivalentes suivantes
 est vérifiées :
 \begin{itemize}
 
 \item $\Delta^\mu$ n'a pas de fonction de Green positive
 \item Il existe une suite $(u_k)_k$ de fonctions lisses à support compact telle que 
 $$ 0\le u_k\le 1\et \lim_{k\to\infty} u_k=1 \mbox{\ uniformément sur les compacts}$$
 $$ \et
 \lim_{k\to\infty} \int_M |du_k|^2 d\mu=0\  \footnote{Cette hypothèse pouvant être remplacée par
  $\limsup_{k\to\infty} \int_M |du_k|^2 d\mu<\infty $.}$$
 \end{itemize}
 Une lecture attentive de la preuve du lemme précédent montre que si pour la mesure
 $d\mu=|\beta|^2d\vol$, l'opérateur $\Delta^\mu$ est non parabolique alors $\alpha$ est nulle en
 cohomologie $L^2$.

 \subsection{Fin de la preuve des théorèmes A et B}
 Nous pouvons complémenter notre corollaire (\ref{coreu}) et finir la preuve des théorèmes A et B :
 \begin{prop}
 \label{inj3}
 Soit $g_0$ une métrique complète sur $\R ^n$, on suppose que $n>2$ et que pour un point $o$
 fixé de $\R^n$ on a les inégalités de Poincaré :
 $$\forall \varphi\in C^\infty(B(o,r)),\  \int_{B(o,r)} (\varphi-m_r(\varphi))^2d\vol\le
  C r^2 \int_{B(o,r)} |d\varphi|^2d\vol $$
  où on a noté $m_r(\varphi)=\int_{B(o,r)} \varphi d\vol/\vol B(o,r)$ la moyenne de la fonction $\varphi$
  sur la boule $B(o,r)$. On suppose aussi que $(\R^n,g_0)$ vérifie l'inégalité de Sobolev :
  $$\forall \varphi\in C_0^\infty(\R^n) \left(\int_{M} (\varphi)^{2\nu/(\nu-2)}d\vol\right)^{1-2/\nu}\le
  C  \int_{M} |d\varphi|^2d\vol$$
  et on suppose enfin que pour une constante $C$ on a pour tout $r>1$
  $$\vol B(o,r)\le C r^\nu$$
Si $(M,g)$ est isométrique à l'infini à  plusieurs copies de $(\R^n,g_0)$  alors
$$H^1_2(M)\simeq H^1_c(M) \et H^{n-1}_2(M)\simeq H^{n-1}(M).$$
 \end{prop}
 Nous commençcons par montrer le lemme suivant qui est nécéssaire pour démontrer cette proposition :
 \begin{lem} Soit $g_0$ une métrique complète sur $\R ^n$ vérifiant les hypothèses de la proposition
 précédente, si $u$ est une fonction lisse sur $\R ^n$ telle que
 $$\int_{\R^n} |du|_{g_0}^2d\vol_{g_0}<\infty$$ alors il y a une constante $c\in \R$ telle que
 pour tout $r>1$ on ait
 $$\int_{B(o,r)} (u-c)^2d\vol\le C r^2 $$
 \end{lem}
 \proof
 Notons en effet $V(r)=\vol B(o,r)$ et $B_k=B(o,2^k)$. Posons
   $$c_k=\frac{1}{V(2^k)}\int_{B_k} u $$
 on va montrer que la suite $c_k$ converge.
 On a 
 \begin{eqnarray*}
 |c_k-c_{k+1}|&=\frac{1}{V(2^k) V(2^{k+1})}
            \left| \int_{B_k\times B_{k+1}}  \left(u(x)-u(y)\right) dxdy \right|\\
	    &\le \frac{1}{V(2^k) V(2^{k+1})}
            \int_{B_{k+1}\times B_{k+1}} \left |u(x)-u(y)\right| dxdy\\
	&\le \frac{1}{V(2^k)}    \left(\int_{B_{k+1}\times B_{k+1}}  \left|u(x)-u(y)\right|^2
	dxdy\right]^{1/2}\\
	&\le \le \frac{\sqrt{2 V(2^{k+1})}}{V(2^k)}    \left(\int_{B_{k+1}}  \left|u(x)-c_{k+1}\right|^2
	dx\right)^{1/2}\\
	&\le C (2^k)^{1-\nu/2} \\
 \end{eqnarray*}
 Où pour obtenir la dernière inégalité, nous avons utilisé l'inégalité de Poincaré et le fait que
 l'inégalité de Sobolev implique que le volume des boules géodésiques 
 de rayon $r$ est uniformément minorée par $C^{te} r^\nu$ (\cite{carsmf}). 
 Grâce à ceci nous savons que la suite $(c_k)_k$ converge vers $c_\infty$ et
 $|c_k-c_\infty |\le C \varepsilon^k$ où
 $\varepsilon=2^{1-\nu/2}\in ]0,1[$.
 Ainsi nous obtenons
 $$\int_{B(o,2^k)} (u-c_\infty)^2 d\vol\le 2 \int_{B(o,2^k)} (u-c_k)^2 d\vol+2V(2^k) |c_k-c_\infty |^2\le C (2^k)^2$$\endproof
 
{\it Démonstration de la proposition :\ }
 Suivant la proposition (2.4) et le théorème (3.3) de \cite{carduke}, on sait que $(M,g)$ vérifie la même inégalité de Sobolev (avec des constantes différentes) et aussi que
 $$\{0\}\lra H^1_c(M)\lra H^1_2(M).$$
 Soit donc $c\in H^1_2(M)$, on veut montrer que $c$ contient une forme lisse à support compact. 
 On sait par hypothèse qu'il y a un compact $K$ de $M$ tel que 
 $M\setminus K=\cup_{i=1}^b E_i$, où chaque $E_i$ est isométrique à $(\R^n\setminus B(o,R),g_0)$.
 En particulier il y a sur chaque $E_i$ une fonction lisse $f_i\in C^\infty(\R^n\setminus B(o,R))$ telle que
 $\alpha=df_i$. Soit $f$ une fonction lisse sur $M$ qui vaut $f_i$ sur chaque $E_i$ on sait
 que $\alpha-df$ est une $1$ forme à support compact. Il reste donc à montrer que dans la classe de
 cohomologie $L^2$ de $df$ il y a une forme à support compacte.
 Soit donc $u_i$ une fonction lisse sur $\R^n$ qui vaut $f_i$ sur $\R^n\setminus B(o,R)$, grâce au lemme
 précédent nous obtenons l'existence d'une constante $c_i$ telle que
 $$\int_{B(o,r)} (u_i-c_i)^2d\vol\le C r^2 $$
 Grâce à cette estimée et au résultat de J. Jost et K. Zuo (cf. \cite{JZ} et la remarque \ref{jost}) nous pouvons affirmer que si $h$ est une fonction lisse sur
 $M$ qui vaut $c_i$ sur $E_i$, alors $dh$ et $df$ ont même classe de cohomologie $L^2$.
 \endproof
 
 \begin{rem}
 En fait la seule propriété topologique de $\R^n$ que nous avons utilisé est le fait qu'il est simplement connexe à
 l'infini. Ainsi $\R^n$ peut être ici changé en n'importe quelle variété simplement connexe à l'infini.
 \end{rem}
 Ceci ne nous permet pas de traiter le cas des surfaces mais celui est aisé car 
 l'espace des $1$-formes harmoniques $L^2$ d'une surface riemanienne  est un invariant conforme. Et nous
 avons la dichotomie suivante pour $(S,g)$  une surface riemannienne complète de type fini :
 \begin{itemize}
 \item Soit $\cH^1(S,g)$ est de dimension finie et $(S,g)$ est conformément équivalente à une surface
 compacte privé d'un nombre fini de points et $$\cH^1(S,g)\simeq
 \ima[H^1_c(S)\lra H^1(S)].$$
 \item Soit $\cH^1(S,g)$ est de dimension infinie.
 \end{itemize}

Cette proposition nous permet d'achever la preuve du théoréme A car il est bien connu qu'un groupe de Lie
nilpotent simplement connexe vérifie nos hypothèses (cf. par exemple \cite{CSV,Sa}).

 Dans le cas où le laplacien de Hodge-deRham présente un trou spectral en degré $0$ et $1$ nous pouvons
 aussi conclure :
 \begin{prop}Soit $g_0$ une métrique complète sur $\R ^n$ où $n>2$ et on suppose que
$(\R^n,g_0)$ vérifie les inégalités de Poincaré :
 
  $$\forall \varphi\in C_0^\infty(\R^n) \int_{M} |\varphi|^{2}d\vol_{g_0}\le
  C  \int_{M} |d\varphi|_{g_0}^2d\vol_{g_0}$$
  et $$\forall \alpha \in C_0^\infty(T^*\R^n) \int_{M} |\alpha|_{g_0}^{2}d\vol_{g_0}\le
  C  \int_{M} [|d\alpha|_{g_0}^2+ |d^*\alpha|_{g_0}^2 ]d\vol_{g_0}.$$
Si $(M,g)$ est isométrique à l'infini à  plusieurs copies de $(\R^n,g_0)$  alors
$$H^1_2(M)\simeq H^1_c(M) \et H^{n-1}_2(M)\simeq H^1(M).$$
 \end{prop}
\proof Pour cela il suffit de remarquer que la première inégalité de Poincaré implique que
 le volume de $(\R^n,g_0)$ est infini et alors on sait grâce à la proposition (5.1) de \cite{CP} que
 $$\{0\}\lra H^1_c(M)\lra H^1_2(M).$$
 Maintenant, l'hypothèse de trou spectral sur les $1$-formes montrent que si $u$ est une fonction lisse sur $\R
 ^n$ telle que
 $\int_{\R^n} |du|_{g_0}^2d\vol_{g_0}<\infty$, alors il y a une constante $c$ telle que $u-c\in L^2$ ; ce qui
 permet d'adapter (encore plus facilement) la preuve de la proposition (\ref{inj3}).
 \endproof
 \subsection{cohomologie $L^2$ et bouts paraboliques.}

 Le deuxième argument géométrique est le suivant :
 \begin{prop}\label{image} On suppose ici que $(M,g)$ est une variété riemannienne complète isométrique au
 dehors d'un domaine compact $D$ à $$( ]0,\infty[ \times \partial D, dt^2+h_t)$$ où $h_t$ est une
 famille de métrique riemannienne sur $\partial D$ telle que
 \begin{itemize}
 \item si $t>s$ alors $h_t\ge h_s$,
 \item si on note $L(t)=\vol_{n-1}( \{t\}\times \partial D)=\int_{\partial D}\sqrt{\det_{h_0} h_t}
 \,d\vol_{h_0}$ alors $$\int_0^\infty \frac{dt}{L(t)}=\infty.$$
 \end{itemize}
L'image de l'application naturelle $H^k_2(M)\lra H^k(M)$ est alors exactement
 $$\ima [H^k_c(M)\lra H^k(M)].$$
  \end{prop}
 \proof Puisque d'après M. Anderson (\cite{A}), l'espace   $\ima [H^k_c(M)\lra H^k(M)]$ s'injecte dans
 $H^k_2(M)$.
 Il suffit de démontrer que si $\alpha$ est une $k$-forme lisse $L^2$ sur $M$, 	alors tirée en
 arrière sur $\partial D$ elle est nulle en cohomologie. Considèrons donc
 $\beta$ une $(n-1-k)$ forme lisse fermée sur $\partial D$ , on veut montrer que
 $$c:=\int_{\partial D}\beta\wedge i_0^*\alpha=0.$$
 où on a noté $i_t\,:\, \partial D\lra [0,\infty[ \times \partial D$ l'inclusion
 $\theta\mapsto (t,\theta)$.
 On a évidement :
 
 $$|\int_{\partial D}\beta\wedge i_0^*\alpha|=|\int_{\partial D}\beta\wedge i_t^*\alpha|\le
 \|\beta\|_{L^\infty} \int_{\partial D} |i_t^*\alpha|_{h_0} d\vol_{h_0}$$
 Maintenant si on note $J(t, \theta)=\sqrt{\det_{h_0(\theta)} h_t(\theta)}$, comme $h_t\ge h_0$ on a
 $|i_t^*\alpha|_{h_0}(\theta)\le J(t, \theta)|\alpha|_{g}(t,\theta)$, d'où en utilisant l'inégalité de
 Cauchy-Schwarz:
 $$c^2\le \|\beta\|_{L^\infty}^2 \, L(t)\, \int_{\partial D} |\alpha|_{g}^2 d\vol_{h_t}$$
 D'où 
 $$c^2\int_r^R \frac{dt}{L(t)}\le  \|\beta\|_{L^\infty}^2 \int_{[r,R]\times \partial D}
 |\alpha|_{g}^2 d\vol_{g}$$ En faisant tendre $R$ vers l'infini, nous obtenons $c=0$.
 \endproof
 
 Remarquons que cette hypothèse implique que $(M,g)$ est parabolique en particulier en dimension
 $2$, ceci implique que $(M,g)$ est conformément équivalent à une surface compacte $\overline{M}$
 privée d'un nombre fini de points et dans ce cas grâce à l'invariance conforme de l'espace des
 $1$-formes harmoniques $L^2$ on sait que
 $$H^1_2(M)\simeq H^1(\overline{M})\simeq \ima [H^1_c(M)\lra H^1(M)]$$
 Nous voulons maintenant en quelque sorte généraliser ce résultat ; d'après le lemme \ref{inj2}, la nullité de
  $H_2^k(M\setminus D,\partial D) $ et les hypothèses de notre proposition (\ref{image}) suffisent pour savoir que :
 $$H^k_2(M)\simeq \ima [H^k_c(M)\lra H^k(M)]\ .$$

Nous allons montrer :
\begin{thm} Si en plus des hypothèses précédentes la fonction  
$$j(t)=\max_{\theta\in \partial D} J(t,\theta)\int_0^t \frac{1}{J(s,\theta)} ds$$ vérifie 
$$\int_0^\infty \frac{dt}{j(t)}= \infty$$
alors nous avons 
$$H^k_2(M)\simeq \ima [H^k_c(M)\lra H^k(M)]$$
\end{thm}
\proof
Soit donc $c\in H_2^k(M\setminus D,\partial D)$ et $\alpha$ une forme lisse fermée représentant
$c$, on a 
$\alpha=d\beta$ où
$$\beta=\int_{-\infty}^0 (\phi^\tau)^* {\rm int}_{X} \alpha\, d\tau$$
où $X$ est le champ de vecteurs $t \frac{\partial}{\partial t}$ sur $]0,\infty[\times \partial D$
et $\phi^\tau$ son flot : $$\phi^\tau (t,\theta)=(e^{\tau}t,\theta).$$
 Compte tenu que l'on suppose que $t\mapsto h_t $ est une famille
croissante de métrique riemannienne sur $\partial D$ on a
$$|\beta(t,\theta)|_g\le \int_0^t |\alpha|_g(s,\theta)ds$$
Maintenant si $l(t,\theta)=\int_0^t J(s,\theta)^{-1} ds$, on a pour toute fonction $u$ telle que $u'\in
L^2(\R,J(t,\theta)dt)$ et
$u(0)=0$
\begin{equation}
\label{hardy}
\frac14\int_0^\infty u^2(t) \left(\frac{l'(t,\theta)}{l(t,\theta)}\right)^2 J(t,\theta) dt\le \int_0^\infty |u'|(t)^2 (t)J(t,\theta) dt
\end{equation}
On obtient cette inégalité d'abord pour tout les $u\in C^\infty_0(]0,\infty[)$ ; en posant
$u(t)=v(t)\,\sqrt{l(t)}$ et en intégrant par partie on obtient 
$$\int_0^\infty |u'(t)|^2 J(t,\theta) dt
\ge\,\frac14 \int_0^\infty |v(t)|^2 (t)\,\frac{1}{l(t,\theta)  J(t,\theta)}\, dt.$$

En appliquant cette inégalité à $u(t)= \int_0^t |\alpha|_g(s,\theta)ds$ et en intégrant
l'inégalité obtenue sur $\partial D$, nous obtenons que 
$$\int_{M\setminus D} \frac{|\beta|^2}{j(t)^2}d\vol <\infty$$ ce qui permet de conclure grâce au
lemme (\ref{poid}).

\endproof
\subsection{Applications}
Ce dernier résultat nous permet de re-démontrer dans une preuve unifiée de nombreux résultats :
\begin{itemize}
\item Lorsque $(M,g)$ est à bouts cymindriques (i.e. lorsque pour tout $t$ : $h_t=h_0$) alors nous
retrouvons un résultat de M. Atiyah, V. Patodi, I. Singer (\cite{APS}).
\item Lorsque $M$ est difféomorphe à l'intérieur d'une variété $\overline M$ dont le bord est constitué de deux
hypersurfaces $X_1$ et $X_2$ s'intersectant suivant une sous variété de codimension $2$ dans $\overline M$ et
lorsque $g$ est une $b$ métrique associée à cette variété à coins i.e. 
$$g=\left(\frac{dx_1}{x_1}\right)^2+ \left(\frac{dx_2}{x_2}\right)^2+h$$
où $x_1$ (resp. $x_2$) est une fonction définissant $X_1$ (resp. $X_2$) et $h$ une métrique lisse sur
$\overline M$ ; alors nous retrouvons un résultat de W. Müller (\cite{Mu}).
\item Lorsque $\partial D$ est une fibration sur un cercle 
$$F\lra \partial D\stackrel{\pi}{\longrightarrow} \bS^1,$$ 
$$\mbox{ où  } \, h_t=k+t^2 (\pi^*\theta)^2$$
et $\theta$ une $1$-forme jamais nulle sur $\bS^1$ et $k$ une métrique riemannienne sur $\partial D$
\footnote{On peut aussi supposer que $k$ est un $2$ tenseur symétrique sur $\partial D$ définie positif sur
la distribution verticale et de rang $\dim F$.}; i.e. $g$ est du type "fibred boundary" dans la
terminologie de R. Mazzeo et R. Melrose (\cite{MM}). Nous retrouvons ainsi une petite partie d'un résultat important de
T.Hausel, E.Hunsicker et R.Mazzeo (corollaire 10  de \cite{HHM}) à propos de instantons gravitationnels qui sont ALG.
\item Lorsque $(M,g)$ est plate au dehors d'un compact et parabolique (i.e le volume des boules géodésiques
croit au plus quadratiquement) alors nous retrouvons une partie du résultat de \cite{cargafa}.
\end{itemize}
Ce théorème nous permet aussi de démontrer des résultats nouveaux par exemple :
\begin{cor}Soit $E\lra X$ un fibré vectoriel riemannien sur une variété riemannienne compacte où les
fibres ont des plans euclidiens, équipé de sa métrique riemannienne $L$ vérifie :
 $$\cH^k_2(L)\simeq \ima [H^k_c(L)\lra H^k(L)].$$

\end{cor}\
Ceci répond à une question posée dans (\cite{caroma} p. 102).
\begin{cor}Si $(M,g)$ est une variété riemannienne complète isométrique au
 dehors d'un domaine compact $D$ à $$( ]0,\infty[ \times \partial D, dt^2+h+f(t)^2\theta^2)$$ 
 où $h$ est une métrique riemannienne sur $\partial D$, $\theta$ est une $1$-forme lisse sur $\partial D$ et $f\in
 C^\infty(\R_+,]0,\infty[)$ vérifie 
  $$\int_0^\infty \frac{dt}{f(t)}=\infty$$
alors
 $$H^k_2(M)\simeq \ima [H^k_c(M)\lra H^k(M)].$$
\end{cor}

\end{document}